\documentclass[a4paper]{amsart}

\usepackage{amsmath,amsthm,amsfonts,amssymb,amscd}
\input{xy}
\xyoption{all}

\newenvironment{demo}{\begin{proof}[Proof]}{\end{proof}}
\newtheorem{teo}{Theorem}[section]
\newtheorem{pro}[teo]{Proposition}
\newtheorem{lem}[teo]{Lemma}
\newtheorem{cor}[teo]{Corollary}

\newtheorem*{rem}{Remark}
\newtheorem*{exa}{Example}

\newtheorem*{prob}{Problem}
\newtheorem*{conj}{Conjecture}

\newenvironment{sis}{\left\{\begin{aligned}}{\end{aligned}\right.}

\numberwithin{equation}{section}

\newcommand{\Z}{\mathbb{Z}}
\newcommand{\C}{\mathbb{C}}
\newcommand{\PP}{\mathbb{P}}
\newcommand{\OO}{\mathcal{O}}
\newcommand{\LL}{\mathcal{L}}
\newcommand{\CC}{\mathcal{C}}

\begin{document}

\title{Families of n-gonal curves with maximal variation of moduli}
\author{Sergey Gorchinskiy and Filippo Viviani}
\address{Steklov Mathematical Institute, Gubkina str. 8, 119991 Moscow}
\email{serge.gorchinsky@rambler.ru}
\address{Universita' degli studi di Roma Tor Vergata, Dipartimento di
matematica, via della ricerca scientifica 1, 00133 Rome}
\email{viviani@axp.mat.uniroma2.it}
\thanks{The first author was partially supported by RFFI grants 04-01-00613 and 
05-01-00455.}
\maketitle

\section{Introduction}

In 1963 Manin proved the Mordell conjecture for function fields (see
\cite{Man}):
\emph{Let $K$ be a function field and let $X$ be a nonisotrivial curve of genus
at least $2$ defined
over $K$. Then $X$ has finitely many $K$-rational points.}

Some years later Parshin (in the case of a complete base, \cite{Par}) and
Arakelov
(in the general case, \cite{Ara}) proved the Shafarevich conjecture for function
fields:
\emph{Let $B$ be a nonsingular, projective, complex curve and let $S$ be a
finite subset of points
of $B$. Fix an integer $g\geq 2$. Then there exist only finitely many
nonisotrivial families of
smooth curves of genus $g$ over $B-S$.}
Moreover Parshin showed how Mordell conjecture follows from Shafarevich
conjecture (this is known as Parshin trick).

Two analogous theorems for number fields were proved by Faltings in 1983 (see
\cite{Fal}). Also in this context, the Parshin trick allows to deduce Mordell
conjecture from Shafarevich one. Moreover it allows to give some explicit
estimations on the number of rational points for fixed number field $K$, $g$ and
the set of points of bad reduction $S$ (see \cite{Szp}).

A uniform version of Shafarevich and Mordell conjectures for function fields was
obtained recently
by L. Caporaso (\cite{cap02}). She proved that, in the case of a one-dimensional
base, there is a
uniform bound for the Shafarevich and Mordell conjectures depending only on the
genus of the fiber,
the genus of the base and the cardinality of the set of bad reduction (where the
fiber is not smooth).
She proved analogous uniform results in higher dimension for  "canonically
polarized varieties", that
is smooth varieties $V$ with ample canonical bundle $K_V$. In that case she
found uniform bounds
depending only on the Hilbert polynomial $h(x)$ of the
canonical polarization ($h(n)=\chi(K_V^n)$), the canonical degree of the
subvariety $T$ of bad reduction and the genus $g$ of the fiber.  We mention that
uniform results in the number
field case are still conjectural (the best result in that direction is contained
in \cite{CHM},
where it is shown that this uniformity result would follow from the Lang
conjectures).

In a subsequent paper (\cite{cap03}), L. Caporaso considered smooth irreducible
subvarieties $V$ of
$\PP^r$ obtaining a uniform bound depending on the degree of the subvariety $V$,
on the degree of
locus $T$ of bad reduction and on the genus $g$ of the fiber.  Also she
described an example, due to
J. de Jong,  which shows that in the case where the place $T$ of bad reduction
has codimension $1$
the bound should depend on its degree (while she proved this is not the case if
$T$ has codimension
bigger that $1$). But in the last section of this paper, she considered an
interesting case where one
can obtain a uniform bound independent from the locus of bad reduction, namely
the case of families
with \emph{maximal variation of moduli}.

A family of smooth curves (or more generally of stable curves) of genus
$g$ over a base $V$ is said to have maximal variation of moduli
if the image of the modular map $V\rightarrow M_g$
(or more generally $V\rightarrow \overline{M_g}$) is of maximal
dimension, namely $\mbox{min}\{{\rm dim}(V),3g-3\}$.
This means that the family is a truly varying
family of curves (just the opposite of an isotrivial family where the modular
map is constant and the fibers don't vary at all).

Of particular interest are the families with maximal varation of moduli over
a base of dimension $3g-3$ because then the modular map is
generically finite and dominant. In that case (and for $g\geq 24$, when it is
known
that $M_g$ is of general type) Caporaso
proved that the number of families over a fixed base $V$ as well as the
number of rational sections of every such family is bounded by a constant
that depends only on the base $V$ and on the genus $g$
(see \cite[Prop. 4]{cap03}).

Moreover she proved (\cite[Lemma 5]{cap03}) that if such a family
has a rational section then
the degree of the modular map (which is called the modular degree)
must be a multiple of $2g-2$ and from the proof one deduces that this result
is sharp, namely that there exist such families with modular degree exactly
$2g-2$.

The aim of this paper is to generalize this lemma to families of
$n$-gonal curves, namely curves that have a $n$ to $1$ map to $\PP^1$
(or equivalently a base-point free $g^1_n$).
In order to explain the results we obtained and to give the ideas of our
proofs, we first review the instructive proof of Caporaso's lemma.
\begin{lem}[Caporaso]
Let $V$ be a complex irreducible variety of dimension $3g-3$ and let
$\mathcal{F}
\rightarrow V$
be a family of smooth curves of genus
$g\geq 2$ with maximal variation of moduli. If this family has a rational
section, then the degree of the modular map $V\rightarrow M_g$
(which is generically finite by hypothesis) is a multiple of $2g-2$.
Moreover this result is sharp, i.e. there exist such families with
a section and with modular degree exactly $2g-2$.
\end{lem}
\begin{demo}
Consider the modular map (generically finite by hypothesis)
$\phi_{\mathcal{F}}:V\rightarrow M_g$
associated to our family $f:\mathcal{F}\rightarrow V$.
Restricting to an open subset of $V$ we can assume the map to be finite
and also with the image contained in $M_g^0$, which is the open subset in $M_g$
corresponding to curves without automorphisms. Now suppose that the family has a
section
$\sigma$ (which we can assume to be regular after restricting the base again)
and look at the following diagram:
$$
\xymatrix{ \ar @{} [dr] |{\Box}
\mathcal{F} \ar@<1ex>[d]^f  \ar[r]^{\widetilde{\phi_{\mathcal{F}}}} &
\mathcal{C}^0_g \ar[d]^{\pi}\\
V \ar@<1ex>[u]^{\sigma}\ar[r]_{\phi_{\mathcal{F}}} & M_g^0 }
$$
where $\mathcal{C}^0_g$ is the universal family
over $M_g^0$. If we call $D$ the horizontal divisor on $\mathcal{C}^0_g$
defined by $D:={\rm Im}(\widetilde{\phi_{\mathcal{F}}}\circ \sigma)$, then
the diagram above factors as a composition of two cartesian diagrams as
follows:
$$
\xymatrix{ \ar @{} [dr] |{\Box}
\mathcal{F} \ar@<1ex>[d]^f \ar[r] &\ar @{} [dr] |{\Box}
\mathcal{C}^0_g\times_{M^0_g}D \ar@<1ex>[d]^{\pi'} \ar[r] &
\mathcal{C}^0_g \ar[d]^{\pi} \\
V \ar@<1ex>[u]^{\sigma} \ar[r]_{\widetilde{\phi_{\mathcal{F}}}\circ \sigma} & D
\ar@<1ex>[u]^{\sigma'}
\ar[r]_{\pi_{|D}} & M_g^0}
$$
where $\sigma'$ is the tautological section. This implies that
${\rm deg}(\phi_{\mathcal{F}})$
is divisible by ${\rm deg}(\pi_{|D})$, which is the relative degree of $D$ with
respect to the map
$\pi:\mathcal{C}_g^0\to M_g^0$.
But $D$ defines an element of the relative Picard group of $\mathcal{C}^0_g$
over $M_g^0$ and by the Franchetta conjecture over $\C$ (now a theorem of
Harer \cite{Har} and Arbarello-Cornalba \cite{AC2}),
this relative Picard group is free of rank 1 generated by the relative
dualizing sheaf $\omega_{\mathcal{C}^0_g/M_g^0}$
which has vertical degree $2g-2$. Hence the relative degree of $D$
should be a multiple of $2g-2$ and the same for the degree of the modular
map.\\
Moreover, taking $D$ an effective divisor representing
$\omega_{\mathcal{C}^0_g/M_g^0}$ and
pulling back the universal family above it, we obtain a family with a section
(the
tautological one) and with modular degree exactly $2g-2$.
\end{demo}
So the main ingredients in the proof of this lemma are the existence of
a universal family over $M_g^0$ and the fact (Franchetta's conjecture,
theorem of Harer-Arbarello-Cornalba) that the relative dualizing sheaf
generates the relative Picard group of this family over the base $M_g^0$.
Unfortunately this theorem is known only over the complex numbers because the
proof
relays heavily on the analytic results of Harer (\cite{Har}) and hence the lemma
of Caporaso is
valid only in characteristic $0$.

If one wants to generalize to $n$-gonal curves, one soon realizes that
the \emph{hyperelliptic} case is very different from the higher gonal case
($n\geq 3$).\\
In fact since every hyperelliptic curve has a non-trivial
automorphism, namely the hyperelliptic involution (and for the generic
hyperelliptic curve this is the only non-trivial automorphism),
there doesn't exist a universal family over any open subset of
the moduli space $H_g$ of hyperelliptic curves (see \cite{GV} for
a detailed study of hyperelliptic families and their relation with the coarse
moduli space $H_g$). Thus there is no hope to generalize the method of
the proof of Caporaso's lemma to hyperelliptic curves.
Nevertheless we realized that the problem of the existence of a rational
section for hyperelliptic families is closely related to another important
problem, that is the existence of a global $g^1_2$ for such families, namely
of a line bundle on the family (defined uniquely up to
the pull-back of a line bundle coming from the base) that restricts on every
fiber to the unique $g^1_2$ of the hyperelliptic curve. Although the unicity
of such a $g^1_2$ on a hyperelliptic curve could make one think that
it should extend to a family, actually this is not the case for $g$ odd (in
general),
while it holds for $g$ even!
Moreover in \cite{GV} we proved that the existence of such a $g^1_2$ is
equivalent
to the Zariski local-triviality
of the family of $\PP^1$ for which the initial family of hyperelliptic curves
is a double cover (in fact it's true that every hyperelliptic family is
a double cover of a family of $\PP^1$, or in other words it's true that
the hyperelliptic involution extends to families but it's not true that
it is associated to a global $g^1_2$!).
And then having reduced the problem to the Zariski local-triviality
of this family of $\PP^1$, we use the existence of a universal family of $\PP^1$
over $H_g$ (non locally-trivial!) to deduce a condition on the divisibility
of the modular map. The result we obtain in section $2$ is the following.
\begin{teo}\label{1.2}
Let $V$ be an irreducible variety of dimension $2g-1$ over an algebraically
closed
field of characteristic different from $2$ and let $\mathcal{F}\rightarrow V$
be a family of smooth hyperelliptic curves of genus $g\geq 2$ with maximal
variation of moduli.
If this family has a rational section then the degree
of the modular map $V\rightarrow H_g$ (which is generically finite by
hypothesis) is a multiple of $2$ and this is sharp (namely there exist
such families with modular degree exactly $2$ for any $g$).
\end{teo}
Note that this result is valid in any characteristic (different from $2$, to
avoid
problems in the construction of double covers) and the proof is completely
algebraic.

The situation is quite different in the \emph{higher gonal} case.
In fact for this case it is known (see section $3$) that the generic
$n$-gonal curve (with $n\geq 3$) doesn't have any non-trivial automorphism
and hence over the moduli space $(M_{g,n-gon})^0$ of $n$-gonal curves without
automorphisms there exists a universal family $\mathcal{C}_{g, n-can}$,
namely the restriction of the universal family over $M_g^0$.
Hence to imitate the proof of Caporaso's lemma it remains to determine the
relative Picard group of this universal family over the base.\\
To do this we use a classical construction of Maroni
(see \cite{Mar1}, \cite{Mar2}) which permits to embed a canonical $n$-gonal
curve
inside an $(n-1)$-rational normal scroll. Moreover it is known
(see \cite{Bal}) that for a generic canonical $n$-gonal curve the rational
normal scrolls obtained in this way are all isomorphic and in fact
they are the "generic'' scrolls, namely the ones can specialize to all the
others.
Hence if we fix one of such
scrolls $X$ and consider the Hilbert scheme ${\rm Hilb}^X_{n-can}$
of canonical $n$-gonal curves inside it, then we have a dominant
map of ${\rm Hilb}^X_{n-can}$ to the locus $M_{g,n-can}$ of $n$-gonal curves.
Now using that the fibers of this map are unirational (they are precisely
${\rm Aut}(X)^0$) we prove that the relative Picard groups of the two universal families
$\mathcal{C}^X_{n-can}\rightarrow {\rm Hilb}^X_{n-can}$ and
$\mathcal{C}_{g,n-gon}\rightarrow (M_{g,n-can})^0$ are isomorphic. Moreover it
is possible to deduce from this construction that there is an effective divisor
on $\mathcal{C}_{g,n-gon}$ representing a relative $G^1_n$.
\\
Now observe that on the universal family
$\mathcal{C}^X_{n-can}$ over ${\rm Hilb}^{X}_{n-can}$
there are two natural line bundles induced by cutting with the hyperplane
section $D$ and with the fiber $f$ of the ruling, and this two line bundles
restrict, on each fiber of
the universal family, to the canonical sheaf and the unique $g^1_n$
correspondingly.\\
So we are naturally lead to the following
\begin{conj}[1]
The relative Picard group of $\mathcal{C}^X_{n-can}\rightarrow {\rm Hilb}^X_{n-can}$
is generated by $D$ and $f$.
\end{conj}
From what we said before, this conjecture is equivalent to the following
\begin{conj}[1']
On the universal family $\mathcal{C}_{g,n-gon}$ over $(M_{g,n-gon})^0$ there is
a
line bundle $G^1_n$ (that restricts to the unique $g^1_n$
on the generic $n$-gonal curve) such that the relative Picard group
$\mathcal{R}(\mathcal{C}_{g,n-gon})$
is generated by $G^1_n$ and the relative canonical sheaf $\omega$.
\end{conj}
In particular it would follow from this second conjecture, just imitating the
proof of Caporaso's
lemma, the following weaker
\begin{conj}[2]
Let  $V$ be an irreducible variety of dimension $2g+2n-5$ (with $4\leq 2n-2< g$)
and let $\mathcal{F}\rightarrow V$ be a family of smooth $n$-gonal curves  of
genus $g$ with maximal
variation of moduli. If this family has a rational section then the degree of
the modular map
$V\rightarrow M_{g,n-gon}$ is a multiple of $\gcd\{n,2g-2\}$. Moreover this
number is sharp, namely there is no other natural number $d$ being a nontrivial
multiple of $\gcd\{n,2g-2\}$ such that for any family with maximal variation of
moduli and with a rational section its modular degree should be a multiple of
$d$.
\end{conj}
In the last part of this article, we prove conjecture 1 and hence conjecture 2
in the case
of \emph{trigonal} curves over an arbitrary algebraically closed field.
Unfortunately our argument seems to work only for families of curves lying on a
surface
(as in the trigonal case). We don't know yet how to attack this
problem in general.

\emph {Aknowledgements} This work is the fruit of a cooperation
between the two authors started at the summer school "Pragmatic
2004",  held at the university of Catania from 25 August to 14
September 2004. We thank prof. L. Caporaso who suggested the
problem to us and who followed our progresses providing constant
encouragement. We thank the co-ordinator of Pragmatic prof. A.
Ragusa, the teachers prof. L. Caporaso and prof. O. Debarre, the
collaborators C. Casagrande, G. Pacienza and M. Paun. We thank prof. 
E. Ballico for pointed out the reference \cite{Bal}. The second
author is grateful to A. Rapagnetta for many helpful
conversations.

\section{Families of hyperelliptic curves}

In this section we work over an algebraically closed field of characteristic
different from $2$.
Recall that the moduli scheme $H_g$ parametrizing isomorphism classes of
hyperelliptic curves is an integral subscheme of $M_g$  of dimension $2g-1$ and
can be described as
$$H_g=({\rm Bin}(2,2g+2)-\Delta)/PGL(2)$$
where ${\rm Bin}(2,2g+2)$ is the projective space of binary forms in two variables of
degree $2g+2$,
$\Delta$ is the closed subset over which the discriminant vanishes and $PGL(2)$
acts naturally on
${\rm Bin}(2,2g+2)$ preserving the locus $\Delta$ (see \cite[Chap. IV, Section
1]{GIT}).\\
We want to study families $\mathcal{F}\rightarrow V$ of smooth hyperelliptic
curves of genus $g\geq 2$
(with $V$ irreducible) such that the modular map
$\phi_{\mathcal{F}}:V\rightarrow H_g$ is generically finite, or equivalently
families over a base $V$ of dimension $2g-1$ with maximal variation of moduli.
We want a "sharp" condition on the degree on the modular map assuming the
existence of a rational section for our family.
It turns out that this problem is very closely related to the following
\begin{prob}
Given a family of hyperelliptic curves $\mathcal{F}\rightarrow V$, does there
exist
a line bundle on $\mathcal{F}$ that restricts to the $g^1_2$ of every fiber?
In other words can the $g^1_2$ be defined globally on a family of hyperelliptic
curves?
\end{prob}
The last problem is birational on the base, namely if there exists a
global $g^1_2$ defined on an open subset of the base $V$, then it extends in a 
unique way to a global $g^1_2$ on the whole $V$ (see \cite[prop. 3.4]{GV}).\\
Let's begin with the following
\begin{pro}\label{2.1}
If a family of smooth hyperelliptic curves $\mathcal{F}\rightarrow V$ has a
rational section then it has
also a globally defined $g^1_2$.
\end{pro}
\begin{demo}
Every family of smooth hyperelliptic curves is a $2:1$
cover of a family
$\mathcal{P}$ of $\PP^1$ (see \cite[theo. 3.1]{GV})
$$\xymatrix{
\mathcal{F} \ar[dr]_{\pi}^{2:1} \ar[dd]_f & \\
& \mathcal{P} \ar[dl]^{\PP^1}_p \\
V &
}$$
Now if the family $\mathcal{F}\rightarrow V$ has a rational section $\sigma$
then also the family
$\mathcal{P}\rightarrow V$ has a rational section given by the composition  $\pi
\circ \sigma$. 
This implies that the family of $\PP^1$ is Zariski locally trivial (see
\cite[prop. 2.1]{GV}) and hence it has
a line bundle $\OO_{\mathcal{P}}(1)$ of vertical degree $1$.
Now the line bundle $\pi^*(\OO_{\mathcal{P}}(1))$ is the required globally
defined $g^1_2$.

Explicitly, as a representative of a global $g^1_2$ we can take the divisor 
$\overline{s(V)}+\overline{\sigma(s(V))}$, where $s:V\to \mathcal{F}$ is the 
section, $\sigma$ stays for the involution $\sigma:\mathcal{F}\to\mathcal{F}$ 
corresponding to $\pi:\mathcal{F}\to\mathcal{P}$, and the bar denotes 
Zariski closure. 
\end{demo}
The converse of this is false, namely to have a globally defined $g^1_2$ is
strictly stronger than having
a rational section.  Nevertheless the following is true.
\begin{pro}\label{2.2}
If a family of smooth hyperelliptic curves $\mathcal{F}\rightarrow V$ has a
globally defined $g^1_2$
then, up to restricting to an open subset of the base, we can find another
family
$\mathcal{F'}\rightarrow V$ with the same modular map and admitting a rational
section.
\end{pro}
We will give two proofs of this proposition.
\begin{demo}[I Proof]
In \cite[prop. 3.4]{GV} we proved that the existence of a global $g_2^1$ is equivalent to
the Zariski local triviality
of the  underlying family $\mathcal{P}$ of $\PP^1$ and actually the global
$g^1_2$ is
the pullback of the line bundle $\OO_{\mathcal{P}}(1)$. This shows that  we can
choose an effective
divisor $D$ on $\mathcal{F}$ that represents the global $g^1_2$. Also, taking
$D$ to be the pull-back of a general rational section of $\mathcal{P}\to X$, we
can suppose
$D$ is not entirely contained in the Weierstrass divisor (the ramification
divisor of the map $\mathcal{F}\rightarrow \mathcal{P}$). So, after possibly
restricting to
an open subset of the base $V$,  we can assume that the projection $D\rightarrow
V$ is an $2:1$
\'etale
cover (and let's call $j$ the involution on $D$ that exchanges the two sheets).
Also on the family
$\mathcal{F}$ there is a natural involution $i$ that on every fiber is the
hyperelliptic involution.
Now consider the diagram
$$\xymatrix{ \ar @{} [dr] |{\Box}
D\times_{V} \mathcal{F} \ar[r] \ar@<1ex>[d] \ar@(ul,ur)[]^{j\times
i}&\mathcal{F} \ar[d] \ar@(u,r)[]^i\\
D \ar[r] \ar@/^/[u]^{\sigma'} \ar@(l,d)[]_j& V
}$$
where $\sigma'$ is the tautological section. Since all the involutions commute
with the map, we can
form the quotients obtaining
$$\xymatrix{
& D\times_{V} \mathcal{F}  \ar[dl] \ar[dr] \ar@<1ex>[dd] \ar@(ul,ur)[]^{j\times
i}& \\
\mathcal{F}'=D\times_V \mathcal{F}/_{(i\times j)}\ar@<1ex>[dd] &
\ar @{} [dr] |{\Box} &\mathcal{F} \ar[dd] \ar@(u,r)[]^i\\
&D \ar[dr] \ar[dl] \ar@/^/[uu]^{\sigma'} \ar@(dl,dr)[]_j& \\
V=D/j \ar@/^/[uu]^{\sigma} & & V\\
}$$
where the tautological section $\sigma'$, being compatible with the involutions,
gives rise
to a section $\sigma$ of the new family  $\mathcal{F}'\rightarrow V$ that, by
construction,
has also the same modular map of the original family, q.e.d.
\end{demo}
\begin{demo}[II Proof]
This proof is done by passing to the generic point $\eta=Spec(k(V))$ of $V$.
In \cite[prop. 2.1 and 3.4]{GV} it is shown that the existence of a global $g^1_2$ is equivalent to
the
isomorphism $\mathcal{P}_{\eta}\cong \PP^1_{\eta}$, where $\mathcal{P}$ as
before is the
family of $\PP^1$ underlying $\mathcal{F}$.  In this case we showed also that
$\mathcal{F}_{\eta}$ is a hyperelliptic curve over $k(V)$ whose affine part is
given in $\mathbb{A}^2$
by an equation of the form $a y^2=f(x)$ where $a\in k(V)^*/(k(V)^*)^2$ and
$f(x)$ is an homogeneous
polynomial of degree $2g+2$ whose roots in $\overline{k(V)}$ defines a point in
$H_g$ corresponding
to the image of $\eta$ under the modular map.\\
Now the existence of a rational section of the family $\mathcal{F}\rightarrow V$
is equivalent
to the existence of a rational point of $\mathcal{F}_{\eta}$ over $k(V)$. But
this is achieved very easily
by varying our $a$ without modifying the polynomial $f(x)$ (and thus the modular
map):
for example if we just take $a':=f(x_0)\neq 0$ with $x_0\in k(V)$ then the new
hyperelliptic curve $\mathcal{F}'_{\eta}$
given by the equation $a' y^2=f(x)$ will have an evident rational solution
$(x,y)=(x_0,1)$.
\end{demo}

Thus studying families of hyperelliptic curves having a rational section from 
the point of view of the degree of their modular map is equivalent to the 
studying of families of hyperelliptic curves with a globally defined $g^1_2$. 
The last problem is solved in the last section of \cite{GV} as an application of 
the theory developed there. 


Let us briefly review certain results from \cite{GV} and describe how do they 
provide an answer to the problem in question. 

\underline{RESULTS from \cite{GV}}
\begin{itemize}
\item[(i)] There exists a family $\mathcal{P}_g\rightarrow H_g^0$ of $\PP^1$ and
a horizontal flat and relatively smooth divisor $D_{2g+2}$ of vertical degree
$2g+2$ that is universal,
 in the sense that every other family
$\mathcal{P}\rightarrow V$ of $\PP^1$ endowed with an horizontal flat and
relatively smooth divisor
$D$ of vertical degree $2g+2$ is the pull-back  of this one by a unique map from
the base
$X$ to $H_g$.
In particular given a family $\mathcal{F}\rightarrow V$ of smooth hyperelliptic
curves, the underlying
family $\mathcal{P}$ of $\PP^1$ together with the branch divisor $D$ of the
$2:1$ cover
$\mathcal{F}\rightarrow \mathcal{P}$  fulfills this properties and hence it is
the pull-back
of the couple $(\mathcal{P}_g,D_{2g+2})$ by mean of the modular map
$V\rightarrow H_g^0$.
Moreover the universal $\PP^1$-family $\mathcal{P}_g\rightarrow H^0_g$
is not
Zariski locally trivial (see \cite[theo. 6.5]{GV}).
\item[(ii)] A family $\mathcal{F}\rightarrow V$ of smooth hyperelliptic curves
has a globally
defined $g^1_2$ if and only if the underlying family $\mathcal{P}\rightarrow V$
of $\PP^1$
is Zariski locally trivial (see \cite[prop. 3.4]{GV}).

\item[(iii)] Given a Zariski locally trivial family $\mathcal{P}$ of $\PP^1$ 
over the base $V$, there exists an open subset $U\subset V$ such that 
$\mathcal{P}|_U$ corresponds to a family $\mathcal{F}\rightarrow V$ of 
hyperelliptic curves (see \cite[theo. 3.5]{GV}).

\end{itemize}

For a discussion of the existence of a global $g^1_2$ for families of
hyperelliptic curves, see
also the last section of \cite{MR}. 

Now we can use these results to answer the initial problem.
\begin{teo}\label{2.3}
If a family of smooth hyperelliptic curves $\mathcal{F}\rightarrow V$ with
generically finite modular map
has a globally defined $g^1_2$, then the degree of the modular map is a multiple
of $2$ and this is
sharp in the sense that there exist families with that property.
\end{teo}
\begin{demo}
Note that, after restricting $V$ to an open subset, we can assume that the image
of the modular
map is contained in $H_g^0$. By (ii) above the existence of a global $g^1_2$ is
equivalent to
the Zariski local-triviality of the underlying family of $\PP^1$.
But by (i) the $\PP^1$-family $\mathcal{P}\rightarrow V$ is the 
pull-back via
the modular (finite) map of the universal family $\mathcal{P}_g\rightarrow
H_g^0$, and since the last one is not Zariski locally trivial a necessary 
condition to become trivial is the parity
of the degree of the modular map (see the last section of \cite{GV}).

On the other hand, there exists a map $V\rightarrow H_g^0$ of degree $2$
such that the pull-back of the universal family $\mathcal{P}_g$ is
Zariski locally trivial. Hence by (iii) this map is modular on a suitable open 
subset $U\subset V$, having hence a globally defined $g^1_2$ by (ii), and this 
concludes the proof.
\end{demo}
Now combining theorem \ref{2.3} with propositions \ref{2.1} and \ref{2.2}, we
get  theorem \ref{1.2}.

\section{Families of $n$-gonal curves ($n\geq 3$)}

In this section we study the case of higher $n$-gonal curves.
Before doing this we want to recall (for the convenience of the reader)
some known classical facts on gonal curves.

Let $g\geq 2$ and $n\geq 2$ be integers. Inside the moduli space
$M_g$ of curves of genus $g$ let us denote with
$M_{g,n-gon}$ the subset corresponding to curves carrying a $g^1_n$
(i.e. a linear system of degree $n$ and dimension $1$). It is known that:\\
(1) $M_{g,n-gon}$ is a closed irreducible subvariety
of $M_g$ (see \cite{Fulton}).\\
(2) The dimension of $M_{g,n-gon}$ is $\min\{3g-3,2n+2g-5\}$.
In particular every curve of genus $g$ has a $g^1_n$ for $2n-2\geq g$
(see \cite{Segre}).\\
(3) The generic $n$-gonal curve with $n>2$ doesn't have non-trivial
automorphisms; the generic $2$-gonal curve (i.e. hyperelliptic curve)
has only the hyperelliptic involution as non-trivial automorphism.
Hence for $n\geq 3$ there exists a universal family over
$(M_{g,n-gon})^0$ (simply the restriction of the universal family
over the open subset of  $M_g$ of curves without non-trivial
automorphisms).\\
(4) As for the number of $g^1_n$ carried by a generic $n$-gonal curve,
 we have that:

\begin{itemize}

\item (i) If $2n-2<g$, then a generic $n$-gonal curve has only one $g^1_n$
          (see \cite{AC}).
\item (ii) If $2n-2=g$, then the generic $n$-gonal curve has only a finite
      number of $g^1_n$ and this number is equal to
      $\frac{(2n-2)!}{n!(n-1)!}$ (see \cite[pag. 359]{GH}).
\item (iii) If $2n-2>g$, then the dimension of the space of all the $g^1_n$
      is equal to $2n-2-g$ (Brill-Noether theory).\\

\end{itemize}

As explained in the introduction, we will use a classical construction
(due to Maroni \cite{Mar1}, \cite{Mar2}) that allows to embed
a canonical $n$-gonal curve (after having chosed a base-point-free $g^1_n$)
in a rational normal scroll. The construction
is based on the observation that, by the geometric version of the
Riemann-Roch theorem
(see \cite[pag. 12]{ACGH}), on the canonical curve 
a $g^1_n$ is given by effective divisors of degree $n$ lying on
a $(n-2)$-plane. Thus we obtain a ruling of $(n-2)$-planes
parametrized by $\PP^1$. It's a classical result of B. Segre (see \cite{Segre})
that, since $2n-2<g$, this ruling is made of non-intersecting planes and so
they sweep a non-singular rational normal scroll. Observe also that this
construction is canonical for the generic curve, since it has only one
$g^1_n$. From now on, we always assume that $2n-2<g$.

Thus we have embedded our $n$-gonal canonical curve $C$ inside a non-singular
rational normal scroll $X$ of dimension $n-1$ inside $\PP^{g-1}$ (and hence of
degree $g-n+1$). It's well known that the Chow ring of a
rational normal scroll is generated by the hyperplane section $D$ and 
the fiber $f$ of the ruling (see \cite[Chap. 3, Sect. 3]{Ful}):
\begin{equation}\label{Chow}
CH(X)=\Z[D, f]/(f^2, D^{n-2}\cdot f- (g-n+1)).
\end{equation}
By construction it follows that $D$
cuts on our curve the canonical divisor while $f$ cuts the linear system
$g^1_n$.
\begin{lem}\label{class-C}
In the construction above the curve $C$ inside $X$ is rationally equivalent to
$n\cdot D^{n-2}+(n-2)(n-g+1)D^{n-3}\cdot f$.
\end{lem}
\begin{demo}
Since the Chow ring of $X$ is generated by $D$ and $f$, our curve $C$ is
rationally equivalent to $a D^{n-2} + b D^{n-3}\cdot f$. Intersecting with
the fiber $f$ we get
$$n=f\cdot C=f \cdot (a D^{n-2} + b D^{n-3}\cdot f)= a$$
while intersecting with $D$ we get
$$2g-2=D \cdot C=D\cdot (a D^{n-2} + b D^{n-3}\cdot f)=a(g-n+1)+b$$
since $D^{n-1}={\rm deg}(X) D^{n-2}\cdot f={\rm deg}(X) \{pt\}$ and the degree of $X$ is $g-
n+1$.
Solving the two equations we get the desired result.
\end{demo}
Now observe that a $(n-1)$-rational normal scroll is isomorphic abstractly to
the
projectivization of a vector bundle $E$ of rank $n-1$ over $\PP^1$ and it's well
known that every vector bundle on $\PP^1$ splits as direct sum of line bundles
and, multiplying by a line bundle, we can normalize it as $E=\oplus_{i=1}^{n-
1}\OO(-r_i)$
with $0=r_1\leq \cdots \leq r_{n-1}$. Further a necessary and sufficient
condition for
$\PP(E)$ to be embedded as a non-singular rational normal scroll inside $\PP^{g-
1}$
(of degree $g-n+1$) is that, set $N:=\sum_{i=1}^{n-1}r_i$, it holds $N<g-n+1$
and $N\equiv g \pmod {n-1}$. In fact the embedding is the map associated to
the very ample divisor
$c_1(\OO_{\PP(E)}(1))+\left[\frac{g-N}{n-1}-1\right]f$.

The invariants $0=r_1\leq r_2\leq\cdots \leq r_{n-1}$ we obtain via the Maroni construction
are also related to the dimension of the multiples of the $g^1_n$ we start with (see
\cite[section 2]{Sch}). Precisely, if we put  $\eta:=\frac{g-N}{n-1}-1$, it holds that
$$h^0(C, k g^1_n)=
\begin{sis}
&k+1& \text{ if } 0\leq k< \eta &\\
&(j+1)k+1-j \eta -\sum_{t=1}^j r_t& \text{ if } \eta+r_j\leq k< \eta +r_{j+1} &\text{ for } j=1, \cdots, n-2\\
&nk+1-g& \text{ if } \eta+r_{n-1}\leq k.&
\end{sis}$$

Note that there is a finite number of
isomorphic classes of rational normal scrolls inside $\PP^{g-1}$. Hence,
since the locus of the $n$-gonal curves inside $M_g$ is
irreducible, the rational normal scrolls canonically
associated to the generic $n$-gonal curves will be isomorphic.
The next theorem of Ballico (see \cite{Bal}) says that a general
$n$-gonal canonical curve lies inside the generic rational normal 
scroll. 
\begin{teo}[Ballico]\label{multiples}
Let $C$ be a generic $n$-gonal curve of genus $g$ (with $2n-2 < g$) 
and let $g_n^1$ be the unique linear system of dimension $1$ and 
degree $n$. Then 
$$h^0(C, k g_n^1)=\begin{sis}
&k+1 &\text{ if } &\: \: k<\frac{g}{n-1},\\
&n k-g +1 &\text{ if } &\: \: k\geq \frac{g}{n-1}.
\end{sis}$$
\end{teo}
\begin{cor}
Let $r$ be the integer such that $0\leq r<n-1$ and $r\equiv g \mod n-1$.
The rational normal scroll associated via the Maroni construction
to the generic $n$-gonal curve of genus $g$ is abstractely isomorphic
to $\PP(\OO_{\PP^1}^{n-1-r}\oplus \OO_{\PP^1}(-1)^{r})$.
\end{cor}
Since all the embeddings
inside a projective space are conjugated by a projective automorphism, from
now on we can fix a generic rational normal scroll $X^{n-1}_{g-n+1}\subset \PP^{g-1}$
(of dimension $n-1$ and degree $g-n+1$). We will consider the locally closed
subset of the Hilbert scheme of canonical curves inside $X$ consisting of
the curves rationally equivalent to $n\cdot D^{n-2}+(n-2)(n-g+1) D^{n-3}
\cdot f$ and we will call it ${\rm Hilb}^X_{n-can}$. What the preceding theorem
tells us is that the canonical map from ${\rm Hilb}^X_{n-can}$ to $M_{g,n-gon}$ is
dominant. We want to look more closely to the fibers of this map as well as 
to the variety ${\rm Hilb}^X_{n-can}$.

First we need two results about the canonical class of a rational normal scroll
and its automorphism group. 

\begin{lem}\label{canonical}
The canonical class of a rational normal scroll $X^{n-1}_{g-n+1}\subset \PP^{g-1}$
is $K_X=-(n-1)D+(g-n-1)f$.
\end{lem}
\begin{demo}
Recall that $X$ is isomorphic to $\PP(E)$ (with $E$ a vector bundle over 
$\PP^1$ of rank $n-1$ and $c_1(E)=-N$) embedded via the map associated to the 
very ample divisor $D=c_1(\OO_{\PP(E)}(1))+\left[\frac{g-N}{n-1}-1\right]f$.
On $\PP(E)$ we have the two exact sequences  
(let's denote with $T_{\PP(E)/\PP^1}$  the vertical tangent with respect 
to the fibration $\pi$ over $\PP^1$):
\begin{eqnarray*}
0\rightarrow \OO_{\PP(E)}\rightarrow \pi^{*}E\otimes \OO_{\PP(E)}(1)
\rightarrow T_{\PP(E)/\PP^1}\rightarrow 0
\end{eqnarray*}
and 
\begin{eqnarray*}
0\rightarrow T_{\PP(E)/\PP^1}\rightarrow T_{\PP(E)} \rightarrow
\pi^{*}T_{\PP^1} \rightarrow 0
\end{eqnarray*}
or putting them togheter  
\begin{equation}
0\rightarrow \OO_{\PP(E)}\rightarrow \pi^{*}E\otimes \OO_{\PP(E)}(1)
\rightarrow T_{\PP(E)} \rightarrow \pi^{*}T_{\PP^1} \rightarrow 0.
\end{equation}
Taking the first Chern classes in the last exact sequence, we get
\begin{gather*}
c_1(T_{\PP(E)})=c_1(\pi^*T_{\PP^1})+c_1(\pi^{*}E\otimes \OO_{\PP(E)}(1))=\\
\pi^*(c_1(T_{\PP^1}))+(n-1)c_1(\OO_{\PP(E)}(1))+\pi^*(c_1(E))=(n-1)c_1(\OO_{\PP(E)}(1))
+(2-N)f.
\end{gather*}
Hence $K_X=-(n-1)c_1(\OO_{\PP(E)}(1))+(N-2)f$ and substituting 
$D=c_1(\OO_{\PP(E)}(1))+\left[\frac{g-N}{n-1}-1\right]f$ we get the conlusion.
\end{demo}
\begin{pro}\label{automorphism}
Let $X^{n-1}_{g-n+1}\subset \PP^{g-1}$ be a {\rm generic} rational normal scroll. Then:
\begin{itemize}
\item[(i)] ${\rm Aut}(X)$ has dimension $n^2-2n+3$.
\item[(ii)] ${\rm Aut}(X)$ is connected with the exception of the case when 
$X\cong \PP^1 \times \PP^1$ in which case it has two connected component (according to
whether an automorphism exchanges or not the two components of $\PP^1$).
\item[(iii)] ${\rm Aut}(X)$ is rational.
\end{itemize}
\end{pro}
\begin{demo}
Let's first suppose that 
$X\cong \PP(E)\not\cong \PP^1\times\PP^1$
In this case an automorphism of $\PP(E)$ must respect the
fibration $\pi:\PP(E)\to \PP^1$ (since the only subvarieties of
the scroll isomorphic to $\PP^{n-2}$ are the fibers of the map
$\pi$), so that we have a map ${\rm Aut}(\PP(E))\to {\rm Aut}(\PP^1)$. The
kernel of this map is the group ${\rm Aut}(\PP(E)_{\PP^1})$ of
\emph{vertical} automorphisms, i.e. automorphisms of the scroll
that induce the identity on the base of the fibration $\PP(E)\to
\PP^1$. So we have an exact sequence:
\begin{equation}\label{auto-sequence}
0\to {\rm Aut}(\PP(E)_{\PP^1})\to {\rm Aut}(\PP(E)) \to {\rm Aut}(\PP^1).
\end{equation}

The last map is surjective. In fact given an
automorphism $\phi$ of $\PP^1$, from the fact that $\phi$ doesn't change the 
linear class of divisors on $\PP^1$ and that on $\PP^1$
every vector bundle is split, we have that $\phi^*(E)\cong E$.
Therefore there exists an isomorphism $\PP(\phi^*(E))\cong \PP(E)$ 
commuting with $\phi$ on the base.

Moreover, ${\rm Aut}(\PP(E)_{\PP^1})=\PP({\rm Aut}(E))$. In
general the subgroup of vertical automorphisms of the projectivized bundle 
coming from the automorphism of the vector bundle can be
identified with the subgroup of the vertical automorphisms that
preserve the line bundle $\OO_{\PP(E)}(1)$ (use the fact that
$\pi_*(\OO_{\PP(E)}(-1))=E^*$). However using the fact that the base of the
fibration is $\PP^1$, we can prove that the two groups are the
same. In fact an element $\psi \in {\rm Aut}(\PP(E)_{\PP^1})$ should
preserve the relative canonical sheaf which is
$$K_{\PP(E)/\PP^1}=\pi^*(c_1(E)^{-1})\otimes
\OO_{\PP(E)}(-r)$$ (since the relative tangent is
$T_{\PP(E)/\PP^1}={\rm Hom}(\OO_{\PP(E)}(-1),
\pi^*(E)/\OO_{\PP(E)}(-1))$). Since $\psi$ commutes with $\pi$,
then $\psi^*(\OO_{\PP(E)}(-r))=\OO_{\PP(E)}(-r)$ from which we get
the claim since the Picard group of $\PP^1$ doesn't contain
elements of $r$-torsion.

There is a cohomological interepretation of these ideas. Consider the sheaf of 
non-commutative groups $\underline{SL}(E)$ consisting of the automorphisms of 
$E$ of determinant 1 (indeed, the determinant for automorphisms of vector 
bundles is 
well-defined since it doesn't change under the conjugation). Also there is a 
non-commutative sheaf $\underline{{\rm Aut}}(\PP(E)_{\PP^1})$ and the exact sequence 
in the 
\'etale topology (we really need the \'etale topology since we have to extract 
root of degree $r$ from the determinant):
$$
1\to\mu_{r}\to\underline{SL}(E)\to\underline{{\rm Aut}}(\PP(E)_{\PP^1})\to 1, 
$$
where $\mu_r$ denotes the sheaf of roots of unity of degree $r$. 
From the long cohomological sequence we get the inclusion 
${\rm Coker}(SL(E)\to {\rm Aut}(\PP(E)_{\PP^1}))\subset 
H^1_{\acute e t}(\PP^1,\mu_{r})={\rm Pic}(\PP^1)_r=0$. The last map may be also 
defined as follows: it associates to $h\in {\rm Aut}(\PP(E)_{\PP^1})$ the invertible 
sheaf $\LL\in {\rm Pic}(\PP^1)$ such that 
$h^*\OO_{\PP(E)}(1)\otimes\OO_{\PP(E)}(1)^{-1}\cong\pi^*(\LL)$. As it was 
discussed 
above the sheaf $\OO_{\PP(E)}(r)$ must be preserved by $h$, so 
$\LL\in {\rm Pic}(\PP^1)_r$.

Finally let us mention that we could reformulate this fact in a more explicit 
way by passing from the \'etale topology to the Zariski one, replacing 
$\underline{SL}(E)$ by $\underline{GL}(E)=\underline{{\rm Aut}}(E)$, $\mu_r$ by 
$\OO_{\PP^1}^*$, and 
using the inclusion $ H^1_{\acute e t}(\PP^1,\mu_{r})={\rm Pic}(\PP^1)_r\subset 
{\rm Pic}(\PP^1)=H^1_{Zar}(\PP^1,\OO^*_{\PP^1})$. Namely let's take $h\in 
{\rm Aut}(\PP(E)_{\PP^1})$. 
Consider a Zariski open covering of the base by the open subsets $U_i$ over 
which $E$ 
is trivial, and over which there exist elements $g_i\in {\rm Aut}(E|_{U_i})$ 
coinciding with $h$ on $\PP(E|_{U_i})$. 
Let $A_{ij}\in GL(U_i\cap U_j)$ be the transition functions for $E$. Then on the 
intersection $U_i\cap U_j$ we have the equality 
$$\lambda_{ij}g_i=A_{ij}g_jA_{ij}^{-1}$$
for some $\lambda_{ij}\in\OO^*(U_i\cap U_j)$. As for determinants we get the 
equality $\lambda_{ij}^r\det(g_i)=\det(g_j)$ so the cocyle $\lambda_{ij}^r$ is 
trivial, while the cocycle $\lambda_{ij}$ defines an element in ${\rm Pic}(\PP^1)_r$ 
whose triviality implies the existence of $g\in {\rm Aut}(E)$ whose action on $\PP(E)$ 
coincides with $h$. 

To compute the dimension of ${\rm Aut}(\PP(E))$, note that 
$\PP({\rm Aut}(E))$ is
the open subset of $\PP({\rm End}(E))= \PP(H^0(E\otimes E^*))$ over
which the determinant doesn't vanish and therefore, since 
$E=\PP(\OO_{\PP^1}^{n-1-r}\oplus \OO_{\PP^1}(-1)^{r})$, we have
$$
{\rm dim}{\rm Aut}(\PP(E)_{\PP^1})=h^0(\PP^1, E\otimes E^*)-1=(n-1)^2-1
$$  
from which one get part (i) using the exact sequence \ref{auto-sequence}.

The connectedeness and the rationality of ${\rm Aut}(\PP(E))$ follows from
the exact sequence \ref{auto-sequence} since both the first and the last 
group is connected and rational.

Finally, in the case where $X\cong \PP(E)\cong \PP^1\times \PP^1$, 
an automorphism may also exchange the
two fibrations of $\PP^1$, so that we have an exact sequence:
$$0\to {\rm Aut}(\PP(E)\to \PP^1)\to {\rm Aut}(\PP(E)) \to \mathbb{Z}/2\mathbb{Z}\to 0,$$
and for ${\rm Aut}(\PP(E)\to \PP^1)$ we can repeat the same argument of above
which will conclude our proof. 
\end{demo}
\begin{teo}
For a \emph{generic} rational normal scroll $X^{n-1}\subset \PP^{g-1}$, the
scheme
${\rm Hilb}^X_{n-can}$ is smooth and irreducible of dimension $2g+n^2-2$ and the
natural (dominant) map ${\rm Hilb}^X_{n-can}\rightarrow M_{g,n-gon}$ has generic 
fiber isomorphic
to the algebraic subgroup ${\rm Aut}(X)^0$ (the connected component of unity 
inside ${\rm Aut}(X)$).
\end{teo}
\begin{demo}
From the theory of deformations, we known that the tangent space
to the semiversal space for the embedded deformations of a curve $C$ inside
$X$ has dimension $h^0(C, N_{C/X})$ while the space of obstructions
sits inside $H^1(C, N_{C/X})$. Hence at the point $[C]\in {\rm Hilb}^X_{n-can}$,
it holds:
$$h^0(C, N_{C/X})-h^1(C, N_{C/X})\leq {\rm dim}_{[C]}({\rm Hilb}^X_{n-can})
\leq {\rm dim}T_{[C]}({\rm Hilb}^X_{n-can})\leq h^0(C,N_{C/X}).$$   
We will prove that for every $[C]\in {\rm Hilb}^X_{n-can}$, it holds: 
\begin{gather}
\chi(C, N_{C/X})=2g+n^2-2\label{Euler-normal}\\
h^1(C, N_{C/X})=0 \label{H1-normal}
\end{gather}
from which it follows that ${\rm Hilb}^X_{n-can}$ is smooth of dimension 
$2g+n^2-2$.\\
From the exact sequence 
$$0\to T_C \to T_{X|C}\to N_{C/X}\to 0$$
it follows that
$$\chi(C, N_{C/X})=\chi(C, T_{X|C})-\chi(C, T_C)=\chi(C, T_X{|C})-(3-3g).$$
To compute $\chi(T_{X|C})$, we apply the Riemann-Roch theorem for fiber 
bundles, using the fact that, in the Chow ring of $X$, the class
of $C$ is $n D^{n-2}+(n-2)(n-g+1)D^{n-3}\cdot f$ (by lemma \ref{class-C})
and the canonical class of $X$ is $-(n-1)D+(g-n-1)f$ (by lemma \ref{canonical}):  
\begin{equation}\label{Euler-tangent}\begin{split}
\chi(C, T_{X|C})={\rm deg}(c_1(T_{X|C}))+{\rm rk}(T_{X|C})\cdot(1-g)=\\
=-K_X\cdot C+(n-1)(1-g)=n^2+1-g
\end{split}\end{equation}
from which it follows formula \ref{Euler-normal}.

Now consider the diagram 
$$\xymatrix{
0\ar[r]& T_C \ar[r] \ar[d] & T_{X|C} \ar[r] \ar[d] & N_{C/X} \ar[r] &0 \\
& p^*T_{\PP^1}\ar[r]^{id} \ar[d]& (\pi^*T_{\PP^1})_{|C} \ar[d]& &\\
& 0 & 0 && \\ 
}$$
where $\pi:X \to \PP^1$ indicates the projection of the scroll $X$ onto
$\PP^1$ (or in in other words the map associated to the divisor $f$),
$p:C\to \PP^1$ denotes its restriction to $C$ (which is therefore
the map associated to the unique $g^1_n$ of $C$) and the vertical maps
are the differential maps of these morphisms.

Passing to the cohomological exact sequences, the vanishing of $H^1(C, N_{C/X})$
will follow once we prove that
$H^1(C, T_{X|C})\stackrel{\cong}{\longrightarrow} 
H^1(C, p^*T_{\PP^1})$. We will prove the last isomorphism by 
showing that both these groups have the same dimension 
$g-2n+2$ (since we know that the map between them is surjective).

In fact since $2n-2<g$ and $C$ is a general $n$-gonal curve, from 
theorem \ref{multiples} it follows:
$$h^1(C, p^*(T_{\PP^1})=h^1(C, 2 g^1_n)= h^0(C, 2 g^1_n)-\chi(C, 2g_n^1)=
3-(2n+1-g)=g-2n+2.$$

On the other hand it's easy to see that the automorphisms of $X$ come from projectivity
of $\PP^{g-1}$ and the ones that fix the canonical curve $C$ are a finite number 
because ${\rm Aut}(C)$ is finite. Hence by proposition \ref{automorphism} (since 
$X$ is generic) 
$$ h^0(C, T_{X|C})=h^0(X, T_X)={\rm dim}{\rm Aut}(T_X)=n^2-2n+3,$$
which togheter with formula \ref{Euler-tangent} gives $h^1(C, T_{X|C})=g-2n+2$.

Now two generic curves $C$ and $C'$ in ${\rm Hilb}^X_{n-can}$ are isomorphic if and only 
if there exists a projectivity of $\PP^{g-1}$ sending one into the other.
But clearly, since $X$ is canonically attached to the generic $C$, this projectivity
must stabilize $X$ and hence is an automorphism of $X$ that will 
preserve the rational class of $C$ and hence, in view of proposition
\ref{automorphism}, belongs to ${\rm Aut}(X)^0$.

The connectedeness of ${\rm Hilb}^X_{n-can}$ (and hence its irreducibility because 
of its smoothness) follows from the fact that it has a dominant map into
a connected variety with connected fibers.
\end{demo}
So after this construction, we end up with the following situation
$$\xymatrix{
\mathcal{C}^X_{n-can}\ar[r]\ar[d]& \mathcal{C}_{g,n-gon}\ar[d]\\
{\rm Hilb}^X_{n-can}\ar[r]_{\phi}& M_{g,n-gon} }$$ where
$\mathcal{C}^X_{n-can}\rightarrow {\rm Hilb}^X_{n-can}$ is the universal
family over the Hilbert scheme, $\mathcal{C}_{g,n-gon}\rightarrow
M_{g,n-gon}$ is the tautological family over the locus of
$n$-gonal curves (universal over the open subset of curves without
automorphisms) and we know that the canonical map $\phi$ is
dominant with generic fiber isomorphic to ${\rm Aut}(X)^0$. We want to compute
the relative Picard group (i.e. line bundles of the family modulo
pull-back of line bundles of the base) of these two families of
curves. This is called classically the group of \emph{rationally
determined line bundles} (the terminology is due to the fact that
it doesn't change if we restrict the family to an open subset
of the base, \cite[prop.2.2]{cil1})).\\
In \cite{cil1} there are many properties of this group and in particular
we found there a very usefull result that allows to compare the groups
of rationally determined line bundles for two families between which there
is a correspondence which is ``generically unirational''.\\
More precisely, let $\mathcal{F}=(\mathcal{C},S,p)$ and
$\mathcal{F'}=(\mathcal{C'},S',p')$ be two families of schemes over an
irreducible
and smooth base and let's denote with $\mathcal{R}(\mathcal{F})$ and
$\mathcal{R}(\mathcal{F}')$ the two relative Picard groups.
Let $T\subset S\times S'$ be an algebraic correspondence
between $S$ and $S'$ such that the two projections
$\pi:T\rightarrow S$ and $\pi^{'}:T\rightarrow S'$ are dominant.
For any point $x\in S$ we will denote by $S'_{x}$ the closed
subset of $S'$ associated to $x$ by this correspondence, namely
$S'_{x}=\pi^{'}(\pi^{-1}(x))$. We use analogous notation for $S_{x'}$ if
$x'\in S'$.
\begin{teo}[see \cite{cil1}]
If, for every $x$ ranging in an open dense subset of $S$, $S'_{x}$ is
unirational then there is a natural monomorphism of groups
$\mathcal{R}(\mathcal{F'})\hookrightarrow \mathcal{R}(\mathcal{F})$.
If the same hypothesis is true for $S_{x'}$, then we have a natural isomorphism
$\mathcal{R}(\mathcal{F'})\cong \mathcal{R}(\mathcal{F})$.
\end{teo}
The hypothesis of the theorem are valid 
in our case since $\phi:{\rm Hilb}^X_{n-can}\rightarrow
M_{g,n-gon}$ is dominant between irreducible varieties with generic rational fiber 
${\rm Aut}(X)^0$, 
and it gives the isomorphism between the relative Picard groups of our two
families.

Now observe that on the universal family $\mathcal{C}^X_{n-can}$
over ${\rm Hilb}^{X}_{n-can}$ there are two natural line bundles induced
by cutting with the divisors $D$ and $f$ (we will use the same
letters for this two line bundles). Since over each fiber of the
universal family, $D$ restricts to the canonical sheaf while $f$
restricts to the unique $g^1_n$, these two line bundles correspond
to the relative canonical line bundle (defined everywhere and in a
canonical way) and to a globally defined $g^1_n$ (which is well
defined only as an element of the relative Picard group). Now we
make the following
\begin{conj}[1]
The relative Picard group of $\mathcal{C}^X_{n-can}\rightarrow {\rm Hilb}^X_{n-can}$
is generated by $D$ and $f$.
\end{conj}
From what we said before, this conjecture is equivalent to the following
\begin{conj}[1']
On the universal family $\mathcal{C}_{g,n-gon}$ over $(M_{g,n-gon})^0$ there is
a
line bundle $G^1_n$ (that restricts to the unique $g^1_n$
on the generic $n$-gonal curve) such that the relative Picard group
$\mathcal{R}(\mathcal{C}_{g,n-gon})$
is generated by $G^1_n$ and  the relative canonical sheaf $\omega$.
\end{conj}
\begin{cor}
Any line bundle on the universal family $\mathcal{C}^1_{g,n}$ has relative
degree a multiple of $\gcd\{2g-2,n\}$.
\end{cor}
Now imitating words for words the proof of Caporaso's lemma, one proves that
this conjecture
imply the following answer to the original problem.
\begin{conj}[2]
Let  $V$ be an irreducible variety of dimension $2g+2n-5$ (with $4\leq 2n-2< g$)
and let $\mathcal{F}\rightarrow V$ be a family of smooth $n$-gonal curves  of
genus $g$ with maximal
variation of moduli. If this family has a rational section then the degree of
the modular map
$V\rightarrow M_{g,n-gon}$ is a multiple of $\gcd\{n,2g-2\}$. Moreover this
number is sharp, namely there is no other natural number $d$ being a nontrivial
multiple of $\gcd\{n,2g-2\}$ such that for any family with maximal variation of
moduli and with a rational section its modular degree should be a multiple of
$d$.
\end{conj}
The sharpness may be obtained as follows. On $\mathcal{C}^X_{n-can}$ there
exists an effective divisor representing the sheaf $f$: just fix any fiber on
the scroll $X$ and take its intersection with corresponding curves on $X$. This
divisor on $\mathcal{C}^X_{n-can}$ is ${\rm Aut}(X)^0$-equivariant and so we get an
effective divisor on $\mathcal{C}_{g,n-gon}$ being a section of $G^1_n$. Thus we
get a family with maximal variation of moduli, with a rational section and whose
modular degree is equal to $n$. The sharpness easily follows from this and from
the existence of an effective relative canonical divisor on $\mathcal{C}_{g,n-
gon}$.

Now we are going to prove conjecture (1) for \emph{trigonal} curves an arbitrary
algebraically closed field.
\begin{teo}
Conjecture (1) (and hence (1') and (2)) is true for $n=3$.
\end{teo}

To prove this theorem first we will establish a certain rather general
statement.

Consider a smooth projective surface $S$, and let $L$ be a linear system of
divisors on $S$. Let's say that the system $L$ is {\it rather free} if and only
if for a generic curve $C$ from the linear system $L$ and for each point $x\in
C$ we have an equality
$\dim((L|_C)(-x))=\dim(L|_C)-1$. Here by the restriction $L|_C$ we mean the
image of $L$ under the natural map $H^0(S,\OO_S(C))\to H^0(C,\OO_C(C))$.

\begin{exa}
If $L=H^0(S,\LL)$ where $\LL$ is a very ample sheaf then evidently $L$ is rather
free.
\end{exa}

\begin{exa}
If $H^1(S,\OO_S)=0$, $L=H^0(S,\LL)$ and $(\LL.\omega_S)\le -2$ where $\LL$ is an
invertible sheaf on $S$ and $\omega_S$ denotes a canonical sheaf, then $L$ is
rather free.
\end{exa}
To show this fact first consider an exact sequence of sheaves
$$
0\to\OO_S\to\OO_S(C)\to\OO_C(C)\to 0,
$$
associated to a general curve $C$ from $L$. Together with the condition
$H^1(S,\OO_S)=0$ it leads to the equality $L|_C=H^0(C,\OO_C(C))=H^0(C,\LL|_C)$ 
and
so
$(L|_C)(-x)=H^0(C,\LL|_C(-x))$. Now the property of $L$ to be rather free comes
from the inequality
$\rm{deg}_C(\LL|_C(-x))>2g(C)-2$. Indeed
$$
2g(C)-2=(\omega_X\otimes\LL.\LL)<(\LL.\LL)-1,
$$
and
$$
(\LL.\LL)-1=\rm{deg}_C(\LL|_C(-x))).
$$

Now consider the universal curve $\CC$ inside
$\PP(L)\times S=\PP\times S$. Let's denote by $\pi_1$ and $\pi_2$ two natural
projections to $\PP$ and $S$ respectively.

\begin{teo}\label{ratherfree}
If $L$ is rather free then the map
$$
\pi_2^*\colon {\rm Pic}(S)\to {\rm Pic}(\CC)/\pi_1^*({\rm Pic}(\PP))
$$
is surjective.
\end{teo}
\begin{demo}
{\it Step 1.} Fix a generic point $p\in\PP$ and consider the blow-up
$\sigma(\PP)$ of $\PP$ with the center at $p$. Let $\pi$ be the natural
contraction map from $\sigma(\PP)$ to $\PP$. Then $\pi^*(\CC)\to\CC$ is also a
blow-up with the center at the curve $C$ that is a fiber of $\CC$ over $p$. By a
general result (see \cite{cil1}) the relative Picard group doesn't change so we
may concentrate on the new family $\pi^*(\CC)\to \sigma(\PP)$. Indeed if the
total space of the family is regular then the relative Picard group of a family
is invariant when changing the base by its open subset: the isomorphism between
two relative Picard groups is obtained in one direction just by the restriction
and in the other direction by taking the Zariski closure (here we use the
regularity of the total space).

{\it Step 2.} We have the embedding
$\pi^*(\CC)\subset\sigma(\PP)\times S$. Let $f$ denote the composition of this
embedding with the natural map $\sigma(\PP)\times S\to \PP^{N-1}\times S$ where
$N+1=\dim(L)$ and $\PP^{N-1}$ corresponds to lines in $\PP$ passing through $p$.
Then
$f\colon\pi^*(\CC)\to\PP^{N-1}\times S$ is the
blow-up with the center at the subvariety
$R\subset\PP^{N-1}\times S$ of codimension 2 that may be defined in the
following way: a point
$(l,x)\in\PP^{N-1}\times S$ belongs to $R$ if and only if $x\in C\cap C_q$ where
$C_q$ stays for the fiber of $\CC$ over a generic point $q\in l$. Let's remark
that indeed the set $C\cap C_q$ doesn't depend on $q\in l$ because obviously all
the curves in the pencil $l$ are passing through $(C.C)$ points on $C$ obtained
by intersecting it with one of the element of this pencil.

Hence ${\rm Pic}(\pi^*(\CC))$ is generated by
${\rm Pic}(\PP^{N-1}\times X)={\rm Pic}(\PP^{N-1})\times {\rm Pic}(X)$ and the irreducible
components of $f^{-1}(R)$. Besides, the image of the generator of
${\rm Pic}(\PP^{N-1})$ in ${\rm Pic}(\pi^*(\CC))$ is coming from the base $\sigma(\PP)$.

{\it Step 3.} We claim that $R$ is irreducible. There is an embedding
$R\subset\PP^{N-1}\times C$. Also we may identify
$\PP^{N-1}$ with $\PP(L|_C)$. Hence a fiber of $R$ over a point $x\in C$ is
exactly $\PP(L|_C(-x))$. Since $L$ is rather free, $R$ is just the
projectivization of the vector bundle of rank $N-1$ over $C$ whose fiber over
$x\in C$ is the vector space $(L|_C)(-x)$ (which is of dimension $N-1$ by
hypothesis).

Moreover we may easily compute the class of $f^{-1}(R)$ in the relative Picard
group: on the restriction of the family $\pi^*(\CC)$ over the open subset
$\sigma(\PP)-\PP^{N-1}$ the divisor $f^{-1}(R)$ coincides with the pull-back of
$\PP^{N-1}\times C\subset \PP^{N-1}\times S$ so comes from ${\rm Pic}(S)$.

Thus we get the initial statement.

\end{demo}

\begin{rem}
The proof of this theorem is very much inspired by the proof of Theorem 4.2 from
\cite{cil1}, essentially by the idea of considering pencils of curves inside a
family to obtain the information about the relative Picard group (in \cite{cil1}
this method is said to have been already known to Enriques and Chisini).
\end{rem}

Now let's come back to our initial problem and take $S=X$,
$L=H^0(X,3D+(4-g)f)$, where $X$ is a scroll corresponding to trigonal curves of
genus $g$. Since $\omega_S=-2D+(g-4)f$, $f^2=0$, $(D.f)=1$ and $D^2=g-2$ we have
the equality $(\omega_S.(3D+(4-g)f))=-g-8$.
Also $H^1(X,\OO_X)=0$, and so $L$ is rather free by the example before theorem
\ref{ratherfree}.

Another way to see that $L$ is rather free is to show that
$3D+(4-g)f$ is very ample on $X$. In fact in the case $g$ even,
$X$ will be isomorphic to $\mathbb{F}_0=\PP^1\times \PP^1$ with
$D=C_0+\frac{g-2}{2}f$ (in the notation of \cite[V, section
2]{Har}). Hence $3D+(4-g)f=3C_0+\frac{g-2}{2}f$ is very ample by
\cite[V, cor. 2.18]{Har}. In the case $g$ odd, $X$ will be
isomorphic to $\mathbb{F}_1=\PP(\OO\otimes \OO(-1))$ with
$D=C_0+\frac{g-1}{2}f$. Hence $3D+(4-g)f=3C_0+\frac{g+5}{2}f$ is
very ample by the same corollary.

Next note that ${\rm Hilb}_{n-can}^X$ is a Zariski open subset inside $\PP(L)$ because
the Chow groups of scrolls are discrete: rational equivalency and algebraic
equivalency coincide with each other (in the case of divisors it is again the
reflection of the fact that $H^1(X,\OO_X)=0$). Moreover $\CC^X_{n-can}$ is the
restriction of $\CC$ to ${\rm Hilb}^X_{n-can}$. So applying theorem \ref{ratherfree}
and using the birational invariancy of the relative Picard group we get
conjecture 1 for the trigonal case.

\end{document}